\documentclass[twoside,11pt]{article}
\usepackage{graphicx,amsmath,latexsym,amssymb,amsthm}
\usepackage[colorlinks=black,urlcolor=black]{hyperref}
\font\chuto=cmbx10 at 16pt \font\kamy=lcmssb8
at 9pt \font\kam=lcmss8 at 8pt \font\chudaude=cmcsc10

\date{}

\numberwithin{equation}{section}

\begin{document}
\setlength{\unitlength}{1cm}
\begin{picture}(0,0)
\put(-0.7,1.1){\line(1,0){6.85}} \put(-0.7,1.09){\line(1,0){6.85}}
\put(-0.7,1.08){\line(1,0){6.85}}
\put(-0.7,1.07){\line(1,0){6.85}}
\put(-0.7,1.06){\line(1,0){6.85}}
\put(-0.7,1.01){\line(1,0){6.85}}
\put(-0.7,1.00){\line(1,0){6.85}}
\put(-0.65,0.6){\bf\chudaude \textbf{T}hai \textbf{J}ournal of
\textbf{M}athematics}

\put(-0.7,0.04){\line(1,0){6.85}}
\put(-0.7,0.03){\line(1,0){6.85}}
\put(-0.7,-0.02){\line(1,0){6.85}}
\put(-0.7,-0.03){\line(1,0){6.85}}
\put(-0.7,-0.04){\line(1,0){6.85}}
\put(-0.7,-0.05){\line(1,0){6.85}}
\put(-0.7,-0.06){\line(1,0){6.85}}
\put(-0.65,-0.4){\bf\kam
{\textbf{www.math.science.cmu.ac.th/thaijournal}}}
\put(-0.65,-0.8){\bf\kam Online ISSN 1686-0209}

\end{picture}
\vskip1.5cm


\centerline {\bf \chuto Some classical and recent results }

\vskip.2cm

\centerline {\bf \chuto  concerning renorming theory}


\vskip.8cm \centerline {\kamy Amanollah Assadi$^\dag$  and Hadi Haghshenas$^\ddag,$\footnote{{\tt Corresponding author email:haghshenas60@gmail.com (H. Haghshenas)}\\
\\ {\kamy Copyright \copyright\, 2011 by the Mathematical Association of
Thailand. All rights reserved.}}}

\vskip.5cm

\centerline {$^\dag$Department of Pure Mathematics, Birjand
University, Birjand, Iran} \centerline {e-mail : {\tt
aman-assadi@yahoo.com}}


\centerline {$^\ddag$Department of Pure Mathematics, Birjand
University, Birjand, Iran} \centerline {e-mail : {\tt
haghshenas60@gmail.com}}


\vskip.5cm \hskip-.5cm{\small{\bf Abstract :} The problems
connected with equivalent norms lie at the heart of Banach space
theory. This is a short survey on some recent as well as classical
results and open problems in renormings of Banach spaces.



\hrulefill
\section{Introduction.}
Banach space theory is a classic topic in functional analysis. The
study of the structure of Banach spaces provides a framework for
many branches of mathematics like differential calculus, linear
and nonlinear analysis, abstract analysis, topology, probability,
harmonic analysis, etc. The geometry of Banach spaces plays an
important role in Banach space theory. Since it is easier to do
analysis on a Banach space which has a norm with good geometric
properties than on a general space, we consider in this survey an
area of Banach space theory known as \emph{renorming theory}.
Renorming theory is involved with problems concerning the
construction of \emph{equivalent norms} on a vector space with
 nice geometrical properties of convexity or differentiability.
 An excellent monograph containing the main advances on renorming theory until 1993 is
 \cite{14}.\\\\\\\\
We consider only Banach spaces over the reals. Given a Banach
space $X$ with norm $\|.\|$, we denote by $S(X)$ the unit sphere,
and by $X^*$ the dual space with (original) dual norm $\|.\|^{*}$.
All undefined terms and notation are standard and can be found,
for example, in \cite{3, 15, 19, 21, 23, 45, 49}.
\section{Differentiable norms.}
Differentiability of functions on Banach spaces is a natural
extension of the notion of a directional derivative on
$\Bbb{R}^{n}$. A function $f: X\rightarrow \Bbb{R}$ is said to be
\emph{G\^{a}teaux differentiable} at $x \in X$ if there exists a
functional $g \in X^{*}$ such that
$g(y)=\displaystyle{\lim_{t\rightarrow 0}\frac{f(x+ty)-f(x)}{t}},$
for all $y \in X$. In this case, $g$ is called \emph{G\^{a}teaux
derivative} of $f$. If the limit above exists uniformly for each
$y\in S(X)$, then $f$ is called \emph{Fr\'{e}chet differentiable}
at $x$ with \emph{Fr\'{e}chet derivative} $g$. In this paper, most
of our attention will be concentrated on the differentiability of
the norm.
Two important more strong notions of differentiability are
obtained as uniform versions of both Fr\'{e}chet and G\^{a}teaux
differentiability.
The norm $\|.\|$ on $X$ is \emph{uniformly Fr\'{e}chet
differentiable} if $\displaystyle{\lim_{t\rightarrow
0}\frac{\|x+ty\|-\|x\|}{t}}$ exists uniformly for $(x,y)\in
S(X)\times S(X)$. Also, it is \emph{uniformly G\^{a}teaux
differentiable} if for each $y \in S(X)$,
$\displaystyle{\lim_{t\rightarrow
0}\frac{\|x+ty\|-\|x\|}{t}}$ exists uniformly in $x \in S(X)$.\\
Clearly, Fr\'{e}chet differentiability implies G\^{a}teaux
differentiability, but the converse is true only for
finite-dimensional Banach spaces, in general. As an example, the
mapping $f: L^{1}[0,\pi]\rightarrow \Bbb{R}$ defined by
$f(x)=\int_0^\pi{\rm sin}(x(t))dt$ is every where G\^{a}teaux
differentiable, but nowhere Fr\'{e}chet differentiable
\cite{37}.\\
Around the year of 1940, \u{S}mulyan proved his following
fundamental dual characterization of differentiability of norms,
which is used in many basic renorming results.\\\\
\textbf{Theorem 2.1.} \cite[Ch. VIII]{23} \emph{For each $x \in S(X)$, the following are equivalent:\\
\textbf{(i)} $\|.\|$ is Fr\'{e}chet differentiable at $x$.\\
\textbf{(ii)} For all
$(f_{n})_{n=1}^{\infty},(g_{n})_{n=1}^{\infty}\subseteq S(X^{*})$,
if $\displaystyle{\lim_{n\rightarrow \infty}f_{n}(x)=1}$ and
$\displaystyle{\lim_{n\rightarrow \infty}g_{n}(x)=1}$ then
$\displaystyle{\lim_{n\rightarrow
\infty}\|f_{n}-g_{n}\|^{*}=0}$.\\
\textbf{(iii)} Each $(f_{n})_{n=1}^{\infty}\subseteq S(X^{*})$
with $\displaystyle{\lim_{n\rightarrow \infty}f_{n}(x)=1}$ is
convergent in $S(X^{*})$.}\\\\
As a direct application of \u{S}mulyan's theorem, we have the following
corollary:\\\\
   \textbf{Corollary 2.2.} \cite[Ch. VIII]{23} \emph{If the dual norm of $X^{*}$ is Fr\'{e}chet differentiable then $X$ is
    reflexive.}
    \begin{proof} A celebrated and deep theorem of James state that $X$ is reflexive if and only if each nonzero $f\in X^{*}$ attains its
    norm at some $x\in S(X)$. Let $f\in S(X^{*})$ and choose $(x_{n})_{n=1}^{\infty}\in S(X)$
    such that $f(x_{n})\rightarrow1$. By Theorem 2.1, $\displaystyle{\lim_{n\rightarrow\infty}x_{n}=x\in S(X)}$. Therefore
    $f(x)=f(\displaystyle{\lim_{n\rightarrow\infty}x_{n})= \displaystyle{\lim_{n\rightarrow\infty}f(x_{n})=1=\|f\|^{*}.}}$
    If now $f\in X^{*}$ is non-zero, then $\frac{f}{\|f\|^{*}}\in
    S(X^{*})$ and according to the reasoning above there exists $x\in
    S(X)$ such that $\frac{f}{\|f\|^{*}}(x)=1$.\end{proof}
\textbf{Theorem 2.3.} \emph{The following assertions imply the
reflexivity of $X$:\\
\textbf{(i)} The norm of $X$ is uniformly
Fr\'{e}chet differentiable} \cite[page 434]{21}.\\
\emph{\textbf{(ii)} The third dual norm of $X$ is G\^{a}teaux
differentiable} \cite[page 276]{23}.\\\\
\textbf{Theorem 2.4.} \cite[page 275]{23} \emph{If $X$ is
separable and the second dual norm of $X$ is G\^{a}teaux
differentiable then
$X^{*}$ is separable.}\\\\
Since the separable and reflexive spaces contain numerous nice
structural aspects, they have an important role in our
investigation. In fact, there are many renorming characterizations
of reflexivity and separability, such as follows:\\\\
\textbf{Theorem 2.5.} \cite{26} \emph{If $X$ is reflexive then can
be renormed in such a way that both $X$ and $X^{*}$ have
Fr\'{e}chet differentiable norm.}\\\\
There exist reflexive spaces which do not admit any equivalent
uniformly G\^{a}teaux differentiable norm. This example can be
found in Kutzarova and Troyanski \cite{38}. However, \u{S}mulyan
proved the following positive result. His norm was the predual
norm to the norm defined on $X^{*}$ by
$|\|f\||^{2}=\|f\|^{*}{^{2}}+\sum_{i=1}^{\infty}2^{-i}f^{2}(x_{i}),$ where $(x_{i})_{i=1}^{\infty}$ is dense in $S(X)$.\\\\
\textbf{Theorem 2.6.} \cite[Ch. II]{14} \emph{Any separable space
admits an equivalent uniformly G\^{a}teaux differentiable
norm.}\\\\
\textbf{Theorem 2.7.} (Kadec, Restrepo) \cite{36} \emph{If a
separable space $X$ admits an equivalent Fr\'{e}chet
differentiable norm then $X^{*}$ is separable.}\begin{proof}
Observe that the set $B=\{\|.\|^{'}: x \in X , x\neq0\}$ is
norm-separable, where $\|x\|^{'}$ denotes the derivative of
$\|.\|$ at $x$. The set $B$ contains all norm-attaining
functionals, and is thus norm-dense in $X^{*}$ by the
Bishop-Phelps theorem.\end{proof}
The norm $\|.\|$ on $X$ is called \emph{2-rotund (resp. weakly
2-rotund)} if for every $(x_{n})_{n=1}^{\infty} \subseteq S(X)$
such that $\displaystyle{\lim_{m, n\rightarrow \infty}
\|x_{m}+x_{n}\|}=0,$ there is an $x \in X$ such that
$\displaystyle{\lim_{n\rightarrow\infty}x_{n}=x}$ in the
norm (resp. weak) topology of $X$.\\\\
By using Theorem 2.1, it is proved that if a norm on $X$ is
2-rotund then its dual norm is Fr\'{e}chet differentiable. Also,
if a norm on $X$ is weakly 2-rotund then its dual norm is
G\^{a}teaux differentiable \cite{22}.\\\\
\textbf{Theorem 2.8.} \cite[page 208]{32} \emph{$X$ is reflexive
if and only if it admits an equivalent weakly 2-rotund norm.}\\\\
\textbf{Theorem 2.9.} \cite[page 208]{32} \emph{A separable space
$X$ is reflexive if and only if $X$ admits an equivalent 2-rotund
norm.}\\\\
Note that it is not known if the separability of $X$ has to be
assumed in theorem above.\\\\
The space $X$ is \emph{Hilbert generated space} if there is a
Hilbert space $H$ and a bounded linear operator $T$ from $H$ into
$X$ such that $T(H)$ is dense in $X$.\\\\
\textbf{Theorem 2.10.} \cite{22} \emph{$X$ is a subspace of a
Hilbert generated space if and only if $X$ admits an equivalent
uniformly G\^{a}teaux differentiable norm.}
\section{Asplund spaces.}
It is a well-known theorem that every continuous convex function
on a separable space $X$ is G\^{a}teaux differentiable at the
points of a $G_{\delta}$-dense subset of $X$ \cite[page 384]{21}.
Let $f$ be a continuous convex function on $X$. Then the set $G$
of all points in $X$ where $f$ is Fr\'{e}chet differentiable
(possibly empty) is a $G_{\delta}$ set in $X$ \cite[page 357]{21}.
The space $X$ is said to be \emph{Asplund} if every continuous
convex function on it is Fr\'{e}chet differentiable at each point
of a dense $G_{\delta}$ subset of $X$. There exist many well-known
equivalent characterizations of the Asplund spaces. For example,
$X$ is Asplund if and only if $Y^{*}$ is separable whenever $Y$ is
a separable subspace of $X$.
Every Banach space with a Fr\'{e}chet differentiable norm is
Asplund \cite{40} but, on the other hand, Haydon \cite{34}
constructed Asplund spaces admitting no G\^{a}teaux
differentiable norm.\\\\
\textbf{Theorem 3.1.} \cite{17, 40} \emph{For each separable space
$X$ the following are equivalent:\\
\textbf{(i)} $X^{*}$ is separable.\\
\textbf{(ii)} $X$ is Asplund.\\
\textbf{(iii)} $X$ admits an equivalent Fr\'{e}chet
differentiable norm.\\
\textbf{(iv)} There is no equivalent rough norm on $X$.}\\\\
Recall that the norm $\|.\|$ on $X$ is \emph{rough} if for some
$\varepsilon >0$,$$\displaystyle{\limsup_{h\rightarrow 0}
\frac{1}{\|h\|}(\|x+h\|+\|x-h\|-2)}\geq \varepsilon,$$for every $x
\in S(X)$.\\
By using the canonical norm of $C([0,1])$ as a rough norm, we
obtain that $C([0,1])$ does not admit any Fr\'{e}chet
differentiable norm.
\section{Kadec-Klee property.}
The norm $\|.\|$ on $X$ has \emph{weak-Kadec-Klee} property
provided that whenever $(x_{n})_{n=1}^{\infty}\subseteq X $
converges weakly to some $x \in X$ and
$\displaystyle{\lim_{n\rightarrow \infty}\|x_{n}\|=\|x\|}$, then
$\displaystyle{\lim_{n\rightarrow \infty}\|x_{n}-x\|=0}$. Also, a
dual norm $\|.\|_{*}$ on $X^{*}$ has \emph{$weak^{*}$-Kadec-Klee}
property if $\displaystyle{\lim_{n\rightarrow
\infty}\|f_{n}-f\|_{*}=0}$, whenever
$(f_{n})_{n=1}^{\infty}\subseteq X^{*}$ is weak$^{*}$-convergent
to some $f \in X^{*}$ and $\displaystyle{\lim_{n\rightarrow
\infty}\|f_{n}\|_{*}=\|f\|_{*}.}$\\
The weak-Kadec-Klee norms play an important role in geometric
Banach space theory and its applications.\\\\
\textbf{Theorem 4.1.} \cite[page 422]{21} \emph{\textbf{(i)} Let
$X$ be a separable space. If $X^{*}$ has the weak$^{*}$-Kadec-Klee
property then $X^{*}$ is separable.\\\textbf{(ii)} If $X^{*}$ is
separable then $X$ admits an equivalent norm such that $X^{*}$ has
the weak$^{*}$-Kadec-Klee property.}
\section{Strictly convex spaces.}
One interesting and fruitful line of research, dating from the
early days of Banach space theory, has been to relate analytic
properties of a Banach space to various geometric conditions on
that space. The simplest example of such a condition is that of
strict convexity.\\
The space $X$ (or the norm $\|.\|$ on $X$) is called
\emph{strictly convex} (R) if for $x, y \in S(X)$, $\|x+y\|=2$
implies
$x=y$.\\\\
The following result is a consequence of Theorem 2.1 of
\u{S}mulyan:\\\\
\textbf{Theorem 5.1.} \cite[Ch. VIII]{23} \emph{If a dual norm of
$X^{*}$ is strictly convex (G\^{a}teaux differentiable) then its
predual norm is G\^{a}teaux differentiable (strictly convex).}\\\\
The converse implications in the theorem above are true for
reflexive spaces, but not in general.\\\\
Strict convexity is not preserved by equivalent norms. It is well
known that $\|.\|_{\infty}$ and $\|.\|_{2}$ are equivalent norms
on $\Bbb{R}^{n}$, $\|.\|_{2}$ is strictly convex but
$\|.\|_{\infty}$ is not.\\
A most common strictly convex renorming is based on the following
simple observation. Let $Y$ be a strictly convex space and $T: X
\rightarrow Y$ a linear one-to-one bounded operator; then
$\||x|\|=\|x\|+\|T(x)\|$, $x \in X$, is an equivalent strictly convex norm on $X$.\\\\
\textbf{Theorem 5.2.} \cite{35, 36} \emph{Any  separable space $X$
admits an equivalent norm whose dual norm is strictly convex.}
\begin{proof} Let $\{x_{i}\}_{i=1}^{\infty}$ be dense in $S(X)$. Define a
norm $\||.|\|$ on $X^{*}$ by
$|\|f\||^{2}=\|f\|^{*}{^{2}}+\sum_{i=1}^{\infty}2^{-i}f^{2}(x_{i}).$
It is not hard to show that $|\|.\||$ is a weak$^{*}$ lower
semicontinuous function on $X^{*}$ equivalent with $\|.\|^{*}$.
Hence $|\|.\||$ is the dual of a norm $|.|$ equivalent with
$\|.\|$, and also it is strictly convex.\end{proof}
\section{Locally uniformly convex spaces.}
The concept of a locally uniformly convex norm was introduced by
Lovaglia in \cite{42}.
The space $X$ (or the norm $\|.\|$ on $X$) is said to be
\emph{locally uniformly convex} (LUR) if
$$\displaystyle{\lim_{n\rightarrow \infty} \bigg( 2\|x\|^{2}+
2\|x_{n}\|^{2}-\|x+x_{n}\|^{2} \bigg)}=0
\hspace{0.5cm}\Longrightarrow
\hspace{0.5cm}\displaystyle{\lim_{n\rightarrow
\infty}\|x-x_{n}\|}=0,$$ for any
sequence $(x_{n})_{n=1}^{\infty}$ and $x$ in $X$.\\
Lovaglia showed, as a straightforward consequence of Theorem 2.1,
that the norm of a Banach space is Fr\'{e}chet differentiable if
the dual norm is LUR. The converse does not hold, even up to
renormings. In fact, there exists a space with a Fr\'{e}chet
differentiable norm, which does not admit any equivalent norm with
a strictly convex dual norm \cite{14}. However, in the class of
spaces with unconditional bases, we do have equivalence up to a
renorming.\\Many efforts have been dedicated in the renorming
theory to obtain sufficient conditions for a Banach space to admit
an equivalent LUR norm. In 1979, Troyanski stated the first
characterization of existence of LUR renormings.\\\\
\textbf{Theorem 6.1.}\cite{33} \emph{If $X^{*}$ has a dual LUR
norm then $X$ admits an equivalent LUR
norm.}\\\\
\textbf{Theorem 6.2.} \cite[page 387]{21} \emph{Any space $X$ with
a Fr\'{e}chet differentiable norm which has a G\^{a}teaux
differentiable dual norm admits an equivalent LUR norm.}\\\\
\textbf{Theorem 6.3.} (Kadec) \cite{35} \emph{\textbf{(i)} If $X$
is separable then $X$ admits an equivalent LUR
norm.\\\textbf{(ii)} If $X^{*}$ is separable then $X$ admits an
equivalent norm whose dual norm is LUR.} \begin{proof} We show the
second statement. Suppose that $\{x_{i}\}_{i=1}^{\infty}$ be dense
in $S(X)$ and $\{f_{i}\}_{i=1}^{\infty}$ be dense in $S(X^{*})$.
For $i \in \Bbb{N}$, put $F_{i}=$span$\{f_{1}, f_{2}, ...,
f_{i}\}$. Define a norm $|\|.\||$ on $X^{*}$ by
$$|\|f\||^{2}=\|f\|^{*}{^{2}}+\sum_{i=1}^{\infty}2^{-i}{\rm dist}(f,
F_{i})^{2}+\sum_{i=1}^{\infty}2^{-i}f^{2}(x_{i}).$$It is not hard
to show that $|\|.\||$ is a weak$^{*}$ lower semicontinuous
function on $X^{*}$ equivalent with $\|.\|^{*}$. Hence $|\|.\||$
is the dual of a norm $|.|$ equivalent with $\|.\|$, and also it
is LUR.\end{proof}
The theorem above shows that, in particular, every separable space
admits an equivalent strictly convex norm. By using Theorems 3.1
and 6.3, we see that if $X$ is an Asplund space then $X^{*}$
admits an equivalent LUR norm. The next theorem is a powerful
result of Troyanski
\cite{55}:\\\\
\textbf{Theorem 6.4.} \emph{$X$ admits an equivalent LUR norm if
and only if it admits an equivalent weak-Kadec-Klee norm and an
equivalent strictly convex norm.}\\\\
\textbf{Theorem 6.5.} \cite{51} \emph{Let $Y$ be a closed subspace
of $X$ such that both $Y$ and $X/Y$ admit equivalent norms whose
dual norms are LUR. Then $X$ admits an
equivalent norm whose dual norm is LUR.}\\\\
Let us mention here that the analogue of Theorem 6.5 for
Fr\'{e}chet differentiable norms is still an open question.
Talagrand \cite{53} proved that the corresponding result for
G\^{a}teaux differentiable norms is false.\\\\
Here we offer a characterization of LUR spaces in terms of
Lipschitz separated spaces:\\\\
Given a positive scalar $M$, we will let $L_{X, M}$ be the space
of all functions $f:X \rightarrow \Bbb{R}$ such that $|f(x)-f(y)|
\leq \|x-y\|$ for each $x, y \in X$ and sup$\{|f(x)|: x \in
X\}\leq M$ endowed with the metric $\rho(f,g)=$ sup$\{|f(x)-g(x)|:
x \in X \}$. With this metric, $L_{X, M}$ is a complete metric
space. Given a closed nonempty subset $Y \subseteq X$ and $f \in
L_{X,M}$, we let $L_{f, M}=\{\widetilde{f} \in L_{X, M} :
\widetilde{f}_{|Y}=f\}$. We say $X$ is \emph{Lipschitz separated}
if for every proper closed subspace $Y \subseteq X$ and every $f
\in L_{Y, M}$, we have $\displaystyle{\sup_{\widetilde{f} \in
L_{f, M}} \widetilde{f}(x)}> \displaystyle{\inf_{\widetilde{f} \in
L_{f, M}}
\widetilde{f}(x)}$ for all $x \in X \setminus Y$.\\\\
\textbf{Theorem 6.6.} \cite{6} \emph{The separable space $X$ can
be equivalently renormed so that it is LUR but not Lipschitz
separated.}\\\\
\textbf{Theorem 6.7.} \cite{6} \emph{Any space $X$ with a
separable dual admits an equivalent norm under which $X$ is
Lipschitz separated but not LUR.}
\section{Uniformly convex spaces.}
A Banach space is strictly convex if the midpoint of each chord of
the unit ball lies beneath the surface. In 1936, Clarkson
introduced the stronger notion of uniform convexity. A Banach
space is uniformly convex if the midpoints of all chords of the
unit ball whose lengths are bounded below by a positive number are
uniformly buried beneath the surface. The class of uniformly
convex Banach spaces is very interesting and has numerous
applications.\\\\\\\\The space $X$ (or the norm $\|.\|$ on $X$) is
said to be \emph{uniformly convex} (UR) if for all sequences
$(x_{n})_{n=1}^{\infty}$, $(y_{n})_{n=1}^{\infty} \subseteq X$
$$\displaystyle{\lim_{n\rightarrow \infty}\bigg( 2\|x_{n}\|^{2}+
2\|y_{n}\|^{2}-\|x_{n}+y_{n}\|^{2} \bigg)}=0 \hspace{0.5cm}
\Longrightarrow \hspace{0.5cm} \displaystyle{\lim_{n\rightarrow
\infty}\|x_{n}- y_{n}\|=0}.$$
For example, any Hilbert space is uniformly convex and it can be
shown that $L^{p}$ spaces are uniformly convex whenever
$1<p<\infty$.
We have (UR)$\Rightarrow$(LUR)$\Rightarrow$(R) but the converse is
not true. For example, define a norm $|\|.\||$ on $C([0,1])$ by
$|\|f\||^{2}=\|f\|_{\infty}^{2}+\|f\|_{2}^{2}$, where
$\|.\|_{\infty}$ denotes the standard supremum norm of $C([0,1])$
and $\|.\|_{2}$ denotes the canonical norm of $L^{2}[0,1]$. Then
$|\|.\||$ is strictly convex but not LUR on $C([0,1])$.\\
There is a complete duality between uniform convexity
and uniform Fr\'{e}chet differentiability.\\\\
\textbf{Theorem 7.1.} (Lindenstrauss) \cite{41} \emph{For any
space $X$, the dual norm of $X^{*}$ is uniformly convex if and
only if its predual norm is uniformly Fr\'{e}chet differentiable.
Also, the dual norm of $X^{*}$ is uniformly Fr\'{e}chet
differentiable if and only if its predual norm is
uniformly convex.}\\\\
One of the first theorems to relate the geometry of the norm to
linear topological properties is the following;\\\\
\textbf{Theorem 7.2.} \cite[pages 37-50]{15} \emph{Any uniformly
convex Banach space is reflexive.}
\begin{proof}
Assume that the norm of $X$ is uniformly convex. Then the dual
norm of $X^{*}$ is uniformly Fr\'{e}chet differentiable by Theorem
7.1. Therefore $X^{*}$ is reflexive by Theorem 2.3 and thus X is
reflexive.\end{proof}
The theorem above shows that any Hilbert space is a reflexive
Banach space which is a well-known result in functional analysis.
Note that the class of uniformly convex Banach spaces does not
coincide with the all reflexive Banach spaces: an example of a
reflexive Banach
space which is not uniformly convex can be given.\\\\
Notice that, the space $C([0,1])$ is a separable non reflexive
space. Consequently, $C([0,1])$ admits no equivalent uniformly
convex norm, although, by Theorem 6.3, it does admit an equivalent
LUR norm.\\\\
\textbf{Theorem 7.3.} \cite[Ch. XI]{23} \emph{Any space that
admits an equivalent WUR norm is an Asplund space.}\\\\
The norm $\|.\|$ on $X$ is \emph{weakly uniformly convex} (WUR) if
for all sequences $(x_{n})_{n=1}^{\infty}$,
$(y_{n})_{n=1}^{\infty} \subseteq X$ with
$\displaystyle{\lim_{n\rightarrow \infty} 2\|x_{n}\|^{2}+
2\|y_{n}\|^{2}-\|x_{n}+y_{n}\|^{2} }=0$, then
$\displaystyle{\lim_{n\rightarrow \infty} x_{n}- y_{n}=0}$, in the
weak topology of $X$.\\\\\\\\
Using Theorems 3.1 and 7.3, we have the
following corollary.\\\\
\textbf{Corollary 7.4.} \emph{If the norm of a separable space $X$
is WUR then $X^{*}$ is separable.}
\section{Super-reflexive spaces.}
Given Banach space $Y$, we say that $Y$ is \emph{finitely
representable} in $X$ if for every $\varepsilon > 0$ and for every
finite-dimensional subspace $Z$ of $Y$, there is an isomorphism
$T$ of $Z$ onto $T(Z)\subseteq X$ such that $\|T\| \|T^{-1}\| < 1
+ \varepsilon$. The space $X$ is said to be \emph{super-reflexive}
if every finitely
representable space in $X$ is reflexive.\\
Clearly, every super-reflexive space is reflexive. One of the
well-known super-reflexive spaces are Hilbert spaces. The $L^{p}$
spaces for $1<p<\infty$ are other examples of super-reflexive
spaces. But there are many other super-reflexive
spaces. This class is mapped by the following equivalence.\\\\
\textbf{Theorem 8.1.} \cite[page 436]{21} \emph{The following
assertions are equivalent:\\\textbf{(i)} $X$ is
superreflexive.\\\textbf{(ii)} $X$ admits an equivalent uniformly
convex norm.\\\textbf{(iii)} $X$ admits an equivalent uniformly
Fr\'{e}chet differentiable norm.\\\textbf{(iv)} $X$ admits an
equivalent norm which is uniformly convex and uniformly
Fr\'{e}chet differentiable.}
\section{Mazur intersection property.}
In 1933, it was Mazur who first studied Banach spaces which have
the so called \emph{Mazur intersection property} (MIP): every
bounded closed convex set can be represented as an intersection of
closed balls. A systematic study of this topic was initiated by
Phelps \cite{48}. In 1978, Giles, Gregory and Sims gave some
characterizations of this property \cite{27}. They raised the
question whether every Banach space with the MIP is an Asplund
space. They also characterized the associated property for a dual
space, called the \emph{weak$^{*}$ Mazur intersection property}:
every bounded weak$^{*}$ closed convex set can be represented as an intersection of closed dual balls.\\
Associated with MIP, we have also the following concepts:\\
A set $C$ in $X$ is a \emph{Mazur set} if given $f \in X^{*}$ with
sup $f(C) < \lambda$, then there exists a closed ball $D$ such
that $C \subseteq D$ and sup $f(D) < \lambda$. The space $X$ is
called a \emph{Mazur space} provided that any intersection of
closed balls in $X$
is a Mazur set.\\\\
\textbf{Theorem 9.1.} \cite{31} \emph{Any space with Fr\'{e}chet
differentiable norm satisfies the MIP. Also, every space whose
dual satisfies the MIP is reflexive and each reflexive space with
a Fr\'{e}chet differentiable norm is a Mazur space. Finally, Mazur
spaces with the MIP are Asplund and G\^{a}teaux differentiable.}\\\\
\textbf{Theorem 9.2.} \cite{50} \emph{A Mazur space with the MIP
admits an equivalent Fr\'{e}chet differentiable norm.}\\\\
\textbf{Theorem 9.3.} \cite{51} \emph{Let $Y$ be a subspace of $X$
such that both $Y^{*}$ and $(X/Y)^{*}$ can be renormed to have the
weak$^{*}$ Mazur intersection property. Then $X^{*}$ can be
renormed to have the weak$^{*}$ Mazur intersection property.}
\section{Weakly compactly generated spaces.}
The space $X$ is said to be \emph{weakly compactly generated}
(WCG) if $X$ is the closed linear span of a weakly compact set $K
\subseteq X$. The class of WCG spaces has been intensively studied
during the last forty years and now is in the
core of modern Banach space theory \cite{15, 19, 21}.\\
Recall that the space $X$ is separable if there exists a countable
set $\{x_{n}\}_{n=1}^{\infty}$ with
$\overline{\{x_{n}\}_{n=1}^{\infty}}=X $. An important
characterization of reflexivity is the result that $X$ is reflexive if and only if $B(X)$, the closed unit ball of $X$, is weakly compact.\\
Notice that if $X$ is reflexive, then one may take $K=B(X)$ in the
definition above, whereas if $X$ is separable, with
$\{x_{n}\}_{n=1}^{\infty}$ dense in the $S(X)$, we can take
$K=\{n^{-1} x_{n} \}_{n=1}^{\infty} \bigcup \{0\}$. In this way we
see that both separable and reflexive spaces are WCG.\\\\
\textbf{Theorem 10.1.} \emph{\textbf{(i)} If $X$ is WCG then $X$
admits an equivalent norm that is simultaneously LUR and
G\^{a}teaux differentiable \cite{1, 56}.\\\textbf{(ii)} If $X^{*}$
is WCG then $X$ admits an equivalent norm $|.|$ the dual norm of
which is LUR. In particular, $|.|$ is Fr\'{e}chet differentiable \cite{30}.}\\\\
The first part of theorem above shows that every reflexive space
admits an equivalent norm with weak-Kadec-Klee
property.\\\\
\textbf{Corollary 10.2.} \cite[page 589]{21} \emph{If $X^{*}$ is
WCG then
$X$ is Asplund.}\\\\
The corollary above shows that any reflexive Banach space is
Asplund.\\\\
\textbf{Theorem 10.3.} \cite{1} \emph{If $X$ is WCG then $X^{*}$ admits an equivalent strictly convex dual norm.}\\\\
\textbf{Theorem 10.4.} \cite{5} \emph{If $X$ is WCG and Asplund
then $X^{*}$ admits an equivalent LUR dual norm.}\\\\
If $M$ is a bounded total set in $X$ (i.e., a bounded set $M$ in
$X$ such that $\overline{span}M=X$), we will say that the norm of
$X$ is \emph{dually $M$-2-rotund} if $(f_{n})_{n=1}^{\infty}$ is
convergent to some $f \in B(X^{*})$ uniformly on $M$ whenever
$f_{n} \in S(X^{*})$ are such that $\displaystyle{\lim_{m,
n\rightarrow \infty} \|f_{m}+f_{n}\|}=0$.\\\\
\textbf{Theorem 10.5.} \cite{22} \emph{$X$ is WCG if and only if
$X$ admits an equivalent dually $M$-2-rotund norm for some bounded
total set in $X$.}
\section{Va\v{s}\'{a}k spaces.}
A class of spaces wider than WCG spaces, known as weakly countably
determined or Va\v{s}\'{a}k spaces, was originally defined and
investigated by Va\v{s}\'{a}k.\\
The space $X$ is \emph{Va\v{s}\'{a}k} if there is a sequence
$(B_{n})_{n=1}^{\infty}$ of $weak^{*}$-compact sets in $X^{**}$
such that given $x \in X$ and $u \in X^{**}\backslash X$, there is
$n \in \Bbb{N}$ such that $x \in B_{n}$ and $u \not \in
B_{n}$.\\\\
\textbf{Theorem 11.1.} \cite[Ch. XI]{23} \emph{If $X^{*}$ is
Va\v{s}\'{a}k then $X$ admits
an equivalent Fr\'{e}chet differentiable norm.}\\\\
\textbf{Theorem 11.2.} \cite{43} \emph{Every Va\v{s}\'{a}k space
has an
equivalent norm the dual norm of which is strictly convex.}\\\\
Many of the renorming results for WCG spaces are actually valid
for Va\v{s}\'{a}k spaces. For example, any Va\v{s}\'{a}k space
admits a G\^{a}teaux differentiable norm \cite[Ch. VII]{14}.
Further details can be found in \cite[Ch. VII]{19}.
\section{Uniform Eberlein compact spaces.}
A compact space $K$ is said to be \emph{uniform Eberlein} if $K$
is homeomorphic to a weakly compact subset of a Hilbert
space in its weak topology.\\\\
\textbf{Theorem 12.1.} \cite[page 624]{21} \emph{\textbf{(i)}
$(B(X^{*}),w^{*})$ is uniform Eberlein compact if and only if $X$
admits an equivalent uniformly G\^{a}teaux differentiable
norm.\\\textbf{(ii)} Let $K$ be a compact space. $C(K)$ admits an
equivalent uniformly G\^{a}teaux differentiable norm if and only
if $K$ is uniform Eberlein.}
\section{Bases and renorming theory.}
A \emph{Schauder} basis for $X$ is a sequence
$(x_{n})_{n=1}^{\infty}$ of vectors in $X$ such that every vector
in $X$ has a unique representation of the form
$\displaystyle{\sum_{n=1}^{\infty}a_{n}x_{n}}$ with each $a_{n}$ a
scalar and where the sum is converges in the norm
topology.\\Recall that a series
$\displaystyle{\sum_{n=1}^{\infty}x_{n}}$ is said to be
\emph{unconditionally convergent} if the series
$\displaystyle{\sum_{n=1}^{\infty}x_{n_{i}}}$ converges for every
choice of $n_{1} < n_{2} < n_{3} <\cdots$.\\A Shauder basis
$(x_{n})_{n=1}^{\infty}$ for $X$ is said to be
\emph{unconditional} if for every $x \in X$, its expansion in
terms of the basis $\displaystyle{\sum_{n=1}^{\infty}a_{n}x_{n}}$
converges unconditionally.\\\\
\textbf{Theorem 13.1.} \cite{52} \emph{Let $X$ have an
unconditional basis. Then $X$ admits an equivalent norm with an
LUR dual norm whenever $X$ admits an equivalent Fr\'{e}chet
differentiable
norm.}\\\\
A biorthogonal system $\{x_{i};f_{i}\}_{i \in I}$ in $X \times
X^{*}$ \big(i.e., $f_{i}(x_{j})=\delta_{ij}$ (the Kronecker delta)
for $i, j \in I$\big) is called \emph{fundamental} provided that
$\overline{span}(x_{i})_{i \in I}=X$. A fundamental biorthogonal
system $\{x_{i};f_{i}\}_{i \in I}$ is a \emph{Markushevich} basis
if $(f_{i})_{i \in I}$ separates the points of $X$. A Markushevich
basis $\{x_{i};f_{i}\}_{i \in I}$ is called \emph{shrinking} if
$\overline{span}(f_{i})_{i \in I}=X^{*}$. Clearly, every Schauder
basis is Markushevich. An example of a Markushevich basis that is
not a Schauder basis is the sequence of trigonometric polynomials
$\{e^{i 2 \pi n t}: n=0, \pm1, \pm2, ...\}$ in the
$\widetilde{C}([0,1])$ of complex continuous functions on $[0,1]$
whose values at $0, 1$ are equal, with the sup-norm. If $X^{*}$ is
separable then $X$ has a shrinking Markushevich basis \cite[page
231]{21}.\\
A compact space $K$ is called a \emph{Corson compact space} if $K$
is homeomorphic to a subset $C$ of $[-1,1]^{\Gamma}$, for some set
$\Gamma$, such that each point in $C$ has only a countable number
of nonzero coordinates.\\For example, any metrizable compact is a
Corson compact or any weakly compact set in a Banach space is a
Corson compact or the dual ball for a Va\v{s}\'{a}k space in its
weak$^{*}$-topology is a Corson compact \cite[Ch. VI]{14}.\\A
Banach space $X$ is called \emph{weakly Lindel\"{o}f determined}
(WLD) if $(B(X^{*}),w^{*})$ is a Corson
compact. Every Va\v{s}\'{a}k space is WLD.\\\\
\textbf{Theorem 13.2.} \cite[page 211]{32} \emph{For any space $X$ the following are equivalent:\\
\textbf{(i)} $X$ has a shrinking Markushevich basis.\\
\textbf{(ii)} $X$ is WCG and Asplund.\\
\textbf{(iii)} $X$ is WLD and Asplund.\\
\textbf{(iv)} $X$ is WLD and admits an equivalent norm whose dual
norm is LUR.\\
\textbf{(v)} $X$ is WLD and admits an equivalent Fr\'{e}chet differentiable
norm.}\\\\
\textbf{Theorem 13.3.} \cite{31, 47} \emph{Let $X$ have a
fundamental biorthogonal system $\{x_{i};f_{i}\}_{i \in I}
\subseteq X \times X^{*}$. Then the subspace $Y={\rm
span}(x_{i})_{i \in I}$ admits an
equivalent LUR norm.}\\\\
\textbf{Theorem 13.4.} \cite{47} \emph{Let $X$ have a fundamental
biorthogonal system. Then $X^{*}$ admits an equivalent norm with
the weak$^{*}$-Mazur intersection property (every bounded
weak$^{*}$-closed convex set can be represented as an intersection
of closed dual balls).}\begin{proof} A dual Banach space has the
weak$^{*}$-Mazur intersection property provided its predual has a
dense set of LUR points. Let us consider a biorthogonal system
$\{x_{i};f_{i}\}_{i \in I} \subseteq X \times X^{*}$ such that
$X=\overline{span}(x_{i})_{i \in I}$ and put $Y=$ span$(x_{i})_{i
\in I}$. Using Theorem 13.3, we obtain an equivalent LUR norm
$|.|$ on Y. Let $\|.\|$ be the norm $|.|$ extended to $X$. Then,
the unit ball of $\|.\|$ is the closure of the unit ball of $|.|$.
We claim that $\|.\|$ is LUR at each point of $Y$. Take $y \in Y
\backslash \{0\}$ and a sequence $(x_{n})_{n \in \Bbb{N}}$ in $X$
so that $\displaystyle{\lim_{n\rightarrow \infty} 2\|y\|^{2}+
2\|x_{n}\|^{2}-\|y+x_{n}\|^{2}}=0.$ If we choose $y_{n} \in Y$
with $\|y_{n}-x_{n}\|<\frac{1}{n}$, then
$\displaystyle{\lim_{n\rightarrow \infty} 2|y|^{2}+
2|y_{n}|^{2}-|y+y_{n}|^{2}}=0$ and, hence,
$\displaystyle{\lim_{n\rightarrow \infty}
\|x_{n}-y\|}=\displaystyle{\lim_{n\rightarrow \infty}
|y_{n}-y|}=0$.\end{proof}
In the theorem above, in fact, we prove that every Banach space
with a fundamental biorthogonal system admits an equivalent norm
with a dense set of locally uniformly convex points.\\\\
\textbf{Theorem 13.5.} \cite{51} \emph{Let $X^{*}$ be a dual
Banach space with a fundamental biorthogonal system
$\{x_{i};f_{i}\}_{i \in I} \subseteq X^{*} \times X $. Then $X$
admits an equivalent norm with the MIP.}
\section{Some interesting problems.}
According to the author's knowledge and taste, the following problems in this area arise:\\
(Q1) If the space $X$ has the Radon-Nikodym property (i.e., for
every $\varepsilon>0$ every bounded subset of $X$ has a non-empty
slice of diameter less than $\varepsilon$), does it follow that
$X$ admits an equivalent weak-Kadec-Klee norm? Does it admits an
equivalent strictly convex norm? Is it true that $X$ admits an equivalent LUR norm?\\
(Q2) Does every Asplund space admit an equivalent SSD norm?
\big(If $X$ admits an equivalent SSD norm, then $X$ is Asplund (G.
Godefroy)\big). Recall that the norm $\|.\|$ on $X$ is called
strongly subdifferentiale (SSD) if for each $x \in X$, the
one-sided limit $\displaystyle{\lim_{t\rightarrow
0^{+}}\frac{\|x+ty\|-\|x\|}{t}}$
exists uniformly for $y$ in $S(X)$. Note that the norm $\|.\|$ is Fr\'{e}chet differentiable if and only if it is
G\^{a}teaux differentiable and at the same time SSD.\\
(Q3) Assume that a Banach space $X$ admits an equivalent
G\^{a}teaux differentiable norm and that $X$ admits also an
equivalent SSD norm. Does $X$ admit an equivalent Fr\'{e}chet
differentiable
norm?\\
(Q4) Assume that $X$ is a nonseparable non Asplund space. Does $X$
admit an equivalent norm that is nowhere SSD except at the origin?
For separable non Asplund space the answer is yes.\\ 
(Q5) Assume that the norm of a separable Banach space $X$ has the
property that its restriction to every infinite dimensional closed
subspace $Y \subseteq X$ has a point of Fr\'{e}chet
differentiability on $Y$. Is then $X^{*}$ necessarily
separable?\\
(Q6) Assume that $X$ is Va\v{s}\'{a}k. Does $X$ admit an
equivalent norm that has the following property:
$(f_{n})_{n=1}^{\infty}$ is weak-convergent to some $f \in
B(X^{*})$ whenever $f_{n} \in S(X^{*})$ are such that
$\|f_{n}+f_{m}\|\rightarrow2$ as
$n,m\rightarrow\infty$?\\
(Q7) Assume that $X$ has an unconditional basis and admits a
G\^{a}teaux differentiable norm. Does $X$ admit a norm the dual norm of which is strictly convex?\\\\\\\\
(Q8) (Godefroy) Assume an Asplund space $X$ has a Markushevich
basis $\{x_{i},f_{i}\}_{i \in I}$ with span$\{f_{i}\}_{i \in I}$
norming $X^{*}$. Is $X$ WCG?\\
(Q9) Assume $X$ admits an equivalent Fr\'{e}chet differentiable
norm. Does $X$ admits an equivalent LUR norm?\\
(Q10) It is proved in \cite{6} that every weakly uniformly convex
Banach space is Lipschitz separated. Can a Lipschitz separated
Banach space be equivalently renormed with a weakly uniformly
convex norm? A related question is if Lipschitz separated Banach
space necessarily an Asplund space?\\
(Q11) Is it true that an equivalent Fr\'{e}chet differentiable
norm in a subspace of a separable and reflexive Banach space can
be extended to an equivalent Fr\'{e}chet differentiable norm in
the whole space?\\
(Q12) Assume that for every nonempty closed, bounded and convex
subset $A$ of $X^{*}$ there exists $x \in X$ which attains its
supremum on $A$. Is $X$ Asplund?\\
(Q13) Is it true that an equivalent Fr\'{e}chet differentiable
norm in a subspace of a separable and reflexive Banach space can
be extended to an equivalent Fr\'{e}chet differentiable norm in
the whole space?\\
(Q14) A separable Banach space $X$ is reflexive if and only if $X$
admits an equivalent 2-rotund norm. Is it true in general for
nonseparable spaces?\\
(Q15) Assume that the norm of a separable Banach space $X$ is such
that the restriction of it to every subspace of $X$ is Fr\'{e}chet
differentiable at a point. Must $X^{*}$ be separable?\\
(Q16) Let $X$ be a WLD space and $X$ admits a G\^{a}teaux
differentiable norm.
Does $X$ admit a norm whose dual norm is strictly convex?\\
(Q17) Let $X$ be a WLD space. Is every convex continuous function
on $X$ G\^{a}teaux differentiable at some points?\\
\textbf{Acknowledgement.} We thank the referee for many valuable
remarks which led to significant improvement of the present paper.


\medskip

(Received -- -- --)

(Accepted -- -- --)

\vskip.5cm

\end{document}